\documentclass{ws-ijm}

\usepackage{amssymb}
\usepackage{amsmath}
\usepackage{mathrsfs}                                                           
\usepackage{slashbox}  
\usepackage[all]{xy}
\usepackage{graphicx}
\usepackage{wrapfig} 

\def\b{\bar}

\renewcommand{\t}{\tilde}
\newcommand{\p}{\partial}

\renewcommand{\t}{\tilde}
\newcommand{\rf}[1]{(\ref{#1})}

\begin{document}
\title{HOMOTOPY RELATIONS FOR TOPOLOGICAL VOA}
\author{ANTON M. ZEITLIN}
\address{Department of Mathematics,
Yale University, 442 Dunham Lab, 10 Hillhouse Avenue, New Haven, CT 06511\\
anton.zeitlin@yale.edu, http://math.yale.edu/$\sim$az84
http://www.ipme.ru/zam.html}

\maketitle
\begin{abstract}
We consider a parameter-dependent version of the homotopy
 associative part of the Lian-Zuckerman homotopy algebra and provide
 an interpretation of the multilinear operations of this algebra in terms of 
 integrals over certain polytopes. We explicitly prove the pentagon relation up to
 homotopy and propose a construction of the higher operations.
\end{abstract}



\section{Introduction}

The relations between homotopy algebras and topological vertex operator algebras (TVOA) were investigated by both mathematicians and physicists starting from the early 90s (e.g. \cite{zwiebach}, \cite{lz}). 

One of the most important results in this direction is the one of Lian and Zuckerman, who showed that TVOA properties lead to a construction of the homotopy BV algebra on the space of states of the corresponding TVOA. They started from the bilinear operation, which corresponds to the holomorphic normal ordering of two vertex operators, and then proved that this operation was homotopy commutative and associative. They also defined the bracket operation and proved that it satisfies a homotopy Leibniz and a homotopy Jacobi relations. 
 
It was also conjectured that this structure can be extended to what is known as $BV_{\infty}$-algebra (see e.g. \cite{zuck}, \cite{huang}, \cite{gorbounov},\cite{vallette}). This conjecture was proven for a certain class of TVOA, but so far there has been no explicit construction of the higher operations as it was done for bilinear and trilinear ones in the original article \cite{lz}.

However, it is clear that the Lian-Zuckerman homotopy algebra describes only a part of the structure of the original TVOA, e.g. the 
bilinear operation corresponds only to the normal ordering and doesn't involve other operator product expansion (OPE) coefficients.    
So, it is natural to think of the construction of a homotopy algebra, which embraces all OPE coefficients. This article is an attempt to construct such an algebraic object. Similar ideas were used in \cite{herbst}, \cite{bvym}, \cite{cftym} 
and \cite{khromov}.

First of all, we extend the space of states in TVOA by introducing ''nonlocal'' vertex operators, which correspond 
to OPEs of vertex operators and integrals of OPEs of vertex operators on a positive real line over some polytopes, 
such that they are well-defined under the correlator of the corresponding TVOA. Each of such nonlocal operators will depend on a number of parameters, and each of them will have a parameter which we call $length$. 

It turns out that one can define certain multilinear operations on this space, each of which are associated to certain polytope. For example, the bilinear operation is just an operator product for a given difference between points, the trilinear and quadrilinear ones correspond to integrals over the interval and over the pentagon correspondingly, i.e.  Stasheff polytopes $K_3$ and $K_4$ \cite{stasheff}.

These operations depend on multiple parameters, which correspond to the coordinates on polytopes, and hence are defined on spaces of nonlocal operators of appropriate length, so that the resulting object is well defined under the correlator. 

We explicitly show that bilinear, trilinear and quadrilinear operations we define, together with the differential from TVOA satisfy the relations, which, if we neglect their parameter dependence, coincide with the relations of the $A_{\infty}$ algebra. We conjecture that the $n$-linear operation is related to the integral of a certain operator product over the Stasheff polytope $K_n$ (modulo some corrections which are also integrals over boundary polytopes). 

The structure of the paper is as follows. In Section 2, as a motivation, we recall the Lian-Zuckerman construction of the homotopy associative algebra related to TVOA. In Section 3, we construct the space of nonlocal vertex operators and then introduce the above-mentioned bilinear, trilinear and quadrilinear parameter-dependent operations.  We show that they satisfy the properties of $A_{\infty}$-algebra modulo parameter dependence. In Section 4, we discuss the structure of $n$-linear operations and in the end we indicate a possible relation with the tesselations of a real moduli space.

\section{Lian-Zuckerman homotopy associative algebra}

\noindent{\bf 2.1. Notation and Conventions.} Throughout the paper we will work with vertex operator algebras 
(VOA), using physics notation. Therefore the elements of the VOA's vector space will be referred to as $states$, and $A(z)$ denotes the vertex operator $Y(A,z)$ (see e.g. \cite{benzvi}) corresponding to the state $A$.\\

\noindent{\bf 2.2. Topological VOA and the Lian-Zuckerman homotopy BV algebra.} Topological vertex operator algebra (TVOA) is a vertex superalgebra (see e.g.\cite{benzvi}) that has an additional odd operator 
$Q$ which makes the graded vector space of VOA a chain complex, such that the Virasoro element $L(z)$ is $Q$-exact. The formal definition is as follows (see e.g. \cite{zuck} for more details).\\

\noindent{\bf Definition 2.1.} {\it Let V be a $\mathbb{Z}$-graded vertex operator superalgebra, such that $V=\oplus_i V^i=\oplus_{i,\mu}V^i[\mu]$, where $i$ represents grading of $V$ with respect to conformal weight and $\mu$ 
represents fermionic grading of $V^i$. We call V a topological vertex operator algebra (TVOA) if there exist four elements: $J\in V^1[1]$, $b\in V^2[-1]$, $F\in V^1[0]$, $L\in V^2[0]$, such that 
\begin{eqnarray}
[Q,b(z)]=L(z),\quad  Q^2=0,\quad b_0^2=0,
\end{eqnarray}
where $Q=J_0$, $b(z)=\sum_n b_nz^{-n-2}$, $J(z)=\sum_n J_nz^{-n-1}$, \\
$L(z)=\sum_n L_nz^{-n-2}$, 
$F(z)=\sum_nF_nz^{-n-1}$. 
Here $L(z)$ is the Virasoro element of $V$; the operators $F_0$, $L_0$ are diagonalizable, commute with each other and their egenvalues coincide with fermionic grading and conformal weight correspondingly.}\\

Lian and Zuckerman have observed that each TVOA possesses a rich algebraic structure, namely, the structure of homotopy associative algebra. 
Let us briefly recall their construction. 
One can define the bilinear operation $(\cdot ,\cdot )$, 
which is the cochain map with respect to $Q$: 
\begin{eqnarray}\label{mu}
(A_1,A_2)=Res_z\frac{ A_1(z)A_2}{z}.
\end{eqnarray}
This operation turns out to be homotopy commutative and associative, where homotopies are given by the bilinear and trilinear opertaions below:
\begin{eqnarray}
&&{\rm {m}}(A_1,A_2)=\sum_{i\ge 0}\frac{(-1)^i}{i+1}Res_wRes_{z-w}(z-w)^iw^{-i-1}b_{-1}
(A_1(z-w)A_2)(w)\mathbf{1},
\nonumber\\
&&{\rm{n}}(A_1,A_2,A_3)=\sum_{i\ge 0}\frac{1}{i+1}Res_zRes_w w^iz^{-i-1}(b_{-1}A_1)(z)A_2(w)A_3+\nonumber\\
&&\ \ \ \ \ \ (-1)^{|A_1||A_2|}\sum_{i\ge 0}\frac{1}{i+1}Res_wRes_z z^iw^{-i-1}(b_{-1}A_2)(w)A_1(z)A_3.
\end{eqnarray}
In other words, the following Proposition is true.\\

\noindent {\bf Proposition 2.1.}\cite{lz} {\it The operation $ (\cdot ,\cdot )$ is homotopy commutative and homotopy associative:
\begin{eqnarray}\label{lzrel}
&&Q (A_1,A_2)= (Q A_1,A_2)+(-1)^{|A_1|} (A_1,Q A_2),\nonumber\\
&& (A_1,A_2)-(-1)^{|A_1||A_2|} (A_2,A_1)= \nonumber\\
&&Q{\rm m}(A_1,A_2)+{\rm m}(QA_1,A_2)+(-1)^{|A_1|}{\rm m}(A_1,QA_2),\\
&&\nonumber\\
&& Q{\rm n}(A_1,A_2,A_3)+{\rm n}(QA_1,A_2,A_3)+(-1)^{|A_1|}{\rm n}(A_1,QA_2,A_3)+\nonumber\\
&&(-1)^{|A_1|+|A_2|}{\rm n}(A_1,A_2,QA_3)= ( (A_1,A_2),A_3)- (A_1, (A_2,A_3)).
\end{eqnarray}}

As part of a more general conjecture Lian and Zuckerman suggested that this homotopy associative algebra (homotopy LZ algebra) 
can be extended to the $A_{\infty}$-algebra (see e.g. \cite{stashbook} and Section 4). It was recently proved \cite{gorbounov}, \cite{vallette} for a certain class of TVOAs.  This means that there are "higher homotopies", i.e. 
in the general case there exist nonzero multilinear operations which satisfy the higher associativity relations. However, as can be seen, due to complicated structure of the third order operation it is very hard to prove that conjecture directly.

In the next Section, we consider the modifications of the operations $(,)$ and $n(\cdot,\cdot, \cdot)$ which will depend on the parameters and satisfy the corresponding modified form of relations \rf{lzrel}. We explicitly show how to make one step beyond the construction of Lian and Zuckerman, i.e. we give some explicit construction of the quadrilinear operation. In Section 4, we give a conjecture about the construction of general multilinear operations.

\section{Homotopy associative algebra with parameters from operator products}

\noindent{\bf 3.1. Assumptions and nonlocal vertex operators.} 
In the following we say that some object (some relation) involving vertex operators is defined (true) in the weak sense, when it is defined (true) under the correlator. Moreover, in this section when we write vertex operator $A(t)$, it is assumed that the values of parameter $t$ are real and greater than 0.\\
  
\noindent {\bf 3.2. Nonlocal weakly defined vertex operators on the real line.}  
 In this subsection we introduce a general setup for further constructions, namely we consider extra operators in addition to the ones canonically considered in VOA. 
Suppose we have some vertex algebra $V$. Let us consider the following integral of correlator:
\begin{equation}\label{def1}
\int_D\langle v^*, A_1(t_1+x)A_2(t_2+x)....A_n(t_n+x)v\rangle dt_1\wedge...\wedge dt_n. 
 \end{equation}
Here $v^*\in V^*$ (by $V^*$ we mean restricted dual of $V$), $v\in V$, $A_1,..., A_n\in V$, $x>0$, $D$ is an $n$-dimensional polytope in $\mathbb{R}^n$ which belongs to the region $t_1>t_2>...>t_n\ge0$. This object is well defined for $x>0$, therefore one can consider the nonlocal operator 
\begin{equation}
Y(A_1,..., A_n)_D(x)\equiv \int_DA_1(t_1+x)A_2(t_2+x)....A_n(t_n+x)dt_1\wedge...\wedge dt_n, 
\end{equation}
which is weakly 
defined. 

One can also consider the nonlocal operator 
\begin{equation}\label{def2}
(A_1,A_2)_{\epsilon}(x)=A_1(x+\epsilon)A_2(x), 
 \end{equation}
where $\epsilon>0$. 
Operations \rf{def1}, \rf{def2} can be generalized to the case when operators $A_k$ are 
themselves of the type from \rf{def1}, \rf{def2} with appropriate conditions on the domains of integration. 
For example, if $Y(\t A_1,..., \t A_r)_{\t D}(y)=\int_{D'}A_1(t'_1+y)A_2(t_2+y)....A_n(t'_r+y)dt'_1\wedge...\wedge dt'_r$, 
the expression 
\begin{eqnarray}
Y(A_1,...,A_{k-1}, Y(\t A_1,..., \t A_r)_{\t D},A_{k+1}...A_{n})_D(x)
\end{eqnarray}
is weakly defined if and only if $D$ and $\t D$ are such that $t_{k-1}>\t t_1+t_k$, $\t t_r+t_k>t_{k+1}$ in the domain of integration 
and the expression
\begin{eqnarray}
Y(A_1,...,A_{k-1}, (B,A_k)_{\epsilon},A_{k+1}...A_{n})_D(x)
\end{eqnarray}
is weakly defined if and only  $t_{k-1}>\epsilon+t_k$ in $D$. Similarly, operators
\begin{eqnarray}
&&(Y_D(A_1,...,A_{m}),Y_{D'}(A'_1,...,A'_{r}))_{\alpha}(x),\quad ((A_1,A_2)_{\epsilon},A_3)_{\rho}(x)
\end{eqnarray}
are weakly defined as long as $D'$ is such that $\alpha+\min t_m>\max t'_1> 0$ and $\rho>\epsilon>0$ correspondingly. 

As we see, one can form various nonlocal operators from ordinary vertex operators using operations 
\rf{def1}, \rf{def2}. 
In general, for each nonlocal operator $\mathcal{Y}$ obtained by this procedure, the corresponding correlator $\langle v^*,\mathcal{Y}v\rangle $ (where $v\in V, v^*\in V^*$) can be considered as a polylinear form of the $A_1, A_2, ..., A_n \in V$. In order to be well defined, the arguments of the corresponding vertex operators $A_1(t_1),A_2(t_2),...,A_n(t_n)$ should satisfy the condition $t_1>t_2>...>t_n$. One can define the $length$ of this operator $\mathcal{L}(\mathcal{Y})\equiv \max t_1-\min t_n$ and the $position$ of vertex operator, $\mathcal{P}(\mathcal{Y})=\min t_n$. For example, for the standard vertex operator $A(t)$, $\mathcal{L}(A(t))=0$ and the position is obviously equal to $t$. In general nonlocal operator "fill" the interval $[\mathcal{P}(\mathcal{Y}), \mathcal{P}(\mathcal{Y})+\mathcal{L}(\mathcal{Y})]$.

 In the next section we introduce certain algebraic operations on the space of such nonlocal operators, which will depend on certain polytopes. The domain of these operations will depend on the length of the corresponding operators.\\

\noindent{\bf 3.3. Regularized associativity and commutativity.} 
From now on let $(V,Q)$ be a topological vertex algebra. 
Then we have the following result which is a consequence of the properties of the operator product. \\

\noindent{\bf Proposition 3.1.}{\it 
The operator $Q$ satisfy Leibniz rule with respect to $(\cdot, \cdot)_{\rho}$:
\begin{equation}
Q(A,B)_{\rho}(t)=(QA,B)_{\rho}(t)+(-1)^{|A|}(A, QB)_{\rho}(t),
\end{equation}
where $A,B \in V$. }\\

In this Proposition we assumed that $A,B \in V$. One can see that this statement is true for nonlocal operators too, as long as operation $(\cdot, \cdot)_{\rho}$ is well defined. In all the statements below, for simplicity we assume that the corresponding operators belong to $V$, however one can easily generalize the statements for nonlocal operators of appropriate length and position.\\

\noindent The properties of the OPE and the relation  $[Q,b_{-1}]= L_{-1}$ give the following Lemma.\\

\noindent {\bf Lemma 3.1.}{\it Let $A, B\in V$ and $t>\epsilon>0$. Then the relation below is true in the weak sense:
\begin{eqnarray}
&& (A,B)_{\epsilon}(t)-(-1)^{|B||A|}(B,A)_{-\epsilon}(t)=\nonumber\\
&& Qm_{\epsilon}(A,B)(t)+m_{\epsilon}(QA,B)(t)+(-1)^{|A|}m_{\epsilon}(A,QB)(t),
\end{eqnarray}
where
\begin{eqnarray}
m_{\epsilon}(A,B)(t)=\int^{0}_{-\epsilon}[b_{-1},A(t'+t+\epsilon)B(t'+t)]dt'.
\end{eqnarray}}
In fact, one can show
the 0th mode of $m_{\epsilon}(A,B)(t)$ when expanded in $\epsilon$ coincides with the Lian-Zuckerman operation 
${\rm m}(A,B)$. \\

\noindent {\bf Remark on notation.} We note that the position $\mathcal{P}$ of the operator $m_{\epsilon}(A,B)(t)$ is not $t$, but $t-\epsilon$. However, in the calculations it is useful to use $t$ as a position variable when we act on it by different operations, see e.g. Proposition 3.2. The same will apply to the operator $\t m(t)$, constructed from 
higher commutativity operation, see Lemma 3.2.\\

\noindent The operation $(\cdot, \cdot)_{\epsilon}$ satisfies a "regularized" homotopy associativity relation, see proposition below.\\

\noindent {\bf Proposition 3.2.}{\it Let $A,B,C\in V$. Then the following is true in the weak sense:
\begin{eqnarray}\label{ass}
&&((A,B)_{\alpha_1}C)_{\rho}(t)-(A,(B,C)_{\alpha_2})_{\rho}(t)=\\
&& Qn_{\rho,\alpha_1,\alpha_2}(A,B,C)(t)+n_{\rho,\alpha_1, \alpha_2}(QA,B,C)(t)+\nonumber\\
&& (-1)^{|A|}n_{\rho,\alpha_1, \alpha_2}(A,QB,C)(t)+(-1)^{|A|+|B|}n_{\rho,\alpha_1, \alpha_2}(A,B,QC)(t)\nonumber,
\end{eqnarray}
where $0<\alpha_1,\alpha_2<< \rho$ and 
\begin{eqnarray}
n_{\rho,\alpha_1,\alpha_2}(A,B,C)(t)=n'_{\rho,\alpha_1,\alpha_2}(A,B,C)(t)+(m_{\alpha_1}(A,B),C)_{\rho}(t),
\end{eqnarray}
so that 
\begin{eqnarray}
n'_{\rho,\alpha_1,\alpha_2}(A,B,C)(t)=
\int^{\rho-\alpha_1}_{\alpha_2}(-1)^{|A|}A(t+\rho) [b_{-1},B](t'+t) C(t)dt'.
\end{eqnarray}
}
\begin{proof} To prove this proposition we notice that 
\begin{eqnarray}\label{rel}
&&A(\rho+t)B(\rho-\alpha_1+t)C(t)-A(\rho+t)B(\alpha_2+t)C(t)=\nonumber\\
&&A(t+\rho_2)\int^{\rho-\alpha_1}_{\alpha_2}[L_{-1}, B(t'+t)] dt'C(t).
\end{eqnarray}
Then we use the fact that $[Q, b_{-1}]=L_{-1}$ and apply commutativity relation to the first term in \rf{rel}. As a result, we obtain \rf{ass}.  \end{proof}

Finally in this subsection we prove the following statement, which describes the higher commutativity.\\

\noindent{\bf Lemma 3.2.}{\it The trilinear operator product $n'(A,B,C)$
satisfies the following relation in the weak sense: 
\begin{eqnarray}
&&n'_{\rho, \epsilon_1, \epsilon_2}(A,B,C)(t)+(-1)^{|A||B|+|A||C|+|C||B|}n'_{-\rho, -\epsilon_2, -\epsilon_1}(C,B,A)(t)=\nonumber\\
&&(-1)^{|A||B|}m_{\rho}((B, A)_{-\epsilon_1}, C)(t)-m_{\rho}(A, (B,C)_{\epsilon_2})(t)+\nonumber\\
&&Q\t m_{\rho,\epsilon_1, \epsilon_2}( A, B, C)(t)-\t m_{\rho,\epsilon_1, \epsilon_2}( QA, B, C)(t)-\nonumber\\
&&(-1)^{|A|}\t m_{\rho,\epsilon_1, \epsilon_2}( A, QB, C)(t)-(-1)^{|A|+|B|}\t m_{\rho,\epsilon_1, \epsilon_2}( A, B, QC)(t),
\end{eqnarray}
where $t>\rho>>\epsilon_{1,2}$ and 
\begin{eqnarray}
&&\t m_{\rho,\epsilon_1, \epsilon_2}(A, B, C)(t)= \\
&&[b_{-1}, \int^0_{-\rho}(-1)^{|A|}A(s+\rho+t)\int^{\rho-\epsilon_1}_{\epsilon_2}[b_{-1},B(t'+s+t)]dt'C(t+s)]ds].\nonumber
\end{eqnarray}}
\begin{proof} In order to prove this statement one just need to use the formula below, which is the consequence of simple integration properties
\begin{eqnarray}
&&(-1)^{|A|}A(\rho+t)\int^{\rho-\epsilon_1}_{\epsilon_2}[b_{-1},B(t'+t)]dt'C(t)-\\
&&(-1)^{|A|}A(t)\int^{-\epsilon_1}_{-\rho+\epsilon_2}[b_{-1},B(t'+t)]dt' C(t-\rho)=\nonumber\\
&&[L_{-1}, \int^0_{-\rho}(-1)^{|A|}A(s+\rho+t)\int^{\rho-\epsilon_1}_{\epsilon_2}[b_{-1},B(t'+s+t)]dt'C(t+s)]ds]\nonumber
\end{eqnarray}
and then use the identity $[Q,b_{-1}]=L_{-1}$.\end{proof}

\noindent {\bf 3.4. Pentagon relation.}
Let us consider the folowing nonlocal operator:
\begin{eqnarray}\label{op5}
&&p'_{P}(A_1, A_2, A_3, A_4)(t)=\nonumber\\
&&(-1)^{|A_2|}A_1(\rho+t)\int_{P}[b_{-1},A_2](x+t)[b_{-1},A_3(y+t)]dx\wedge dy A_4(t).
\end{eqnarray}
Here $A_1, A_2, A_3, A_4\in V$, $P=P(\xi, \alpha_1, \alpha_2, \epsilon_1, \epsilon_2)$ is a pentagon in Figure 1, which  is determined by five parameters: $\xi, \alpha_1, \alpha_2, \epsilon_1, \epsilon_2<<\rho$. The relations between parameters are as follows: $\xi, \epsilon_1, \alpha_1<<\epsilon_2, \alpha_2$.

\begin{figure}[hbt] 
%\centering                                                       
\begin{minipage}[c]
{0.9\textwidth}
\centering           
\includegraphics[width=0.7\textwidth]{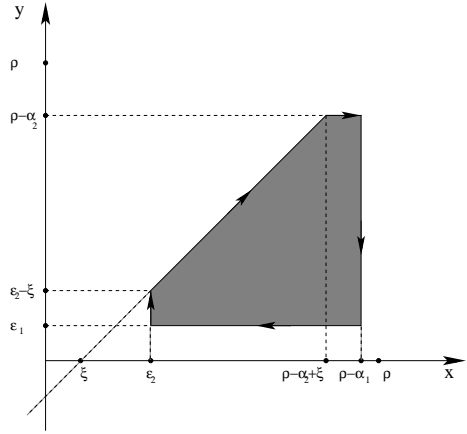}
\end{minipage}%
\begin{minipage}[c]{0.5\textwidth}
\centering
%\begin{tabular}{c}
%hjghjag
%\end{tabular}
\end{minipage}
\caption{Pentagon $P$}                
\end{figure}     

Since 
\begin{eqnarray}
&&[Q,[b_{-1}A_2(x)][b_{-1},A_3(y)]dx\wedge dy]=\nonumber\\
&&\p_xA_2(x)[b_{-1},A_3(y)]dx\wedge dy-(-1)^{|A_2|}[b_{-1}, A_2(x)]\p_yA_3(y)dx\wedge dy=\nonumber\\
&&d(A_2(x)[b_{-1},A_3(y)]dy+(-1)^{|A_2|}[b_{-1}, A_2(x)]A_3(y)dx),
\end{eqnarray}
the following relation is satisfed:
\begin{eqnarray}
&&Qp_{P}'(A_1, A_2, A_3, A_4)(t)-p_{P}'(QA_1, A_2, A_3, A_4)(t)-\nonumber\\
&&(-1)^{|A_1|}p_{P}'(A_1, QA_2, A_3, A_4)(t)-(-1)^{|A_1|+|A_2|}p_{P}'(A_1, A_2, QA_3, A_4)(t)-\nonumber\\
&&(-1)^{|A_1|+|A_2|+|A_3|}p_{P}'(A_1, A_2, A_3, QA_4)(t)=\nonumber\\
&&\nonumber\\
&&(-1)^{|A|_1+|A_2|}A_1(\rho+t)A_2(\epsilon_2+t)\int^{\epsilon_2-\xi}_{\epsilon_1}[b_{-1}, A_3](y+t)dyA_4(t)-\label{1}\\
&&\nonumber\\
&&(-1)^{|A|_1+|A_2|}A_1(\rho+t)A_2(\rho-\alpha_1+t)\int^{\rho-\alpha_2}_{\epsilon_1}[b_{-1},A_3](y+t)dyA_4(t)+\label{2}\\ 
&&\nonumber\\
&&(-1)^{|A_1|}A_1(\rho+t)\int^{\rho-\alpha_1}_{\rho-\alpha_2+\xi}[b_{-1}, A_2](x+t)dxA_3(\rho-\alpha_2)A_4(t)-\label{3}\\ 
&&\nonumber\\
&&(-1)^{|A_1|}A_1(\rho+t)\int^{\rho-\alpha_1}_{\epsilon_2}[b_{-1}, A_2](x+t)dxA_3(\epsilon_1+t)A_4(t)+ \label{4}\\
&&\nonumber\\
&&(-1)^{|A|_1+|A_2|}A_1(\rho+t)\int^{\rho-\alpha_2}_{\epsilon_2-\xi}A_2(y+\xi+t)[b_{-1}, A_3](y+t)dyA_4(t)+\label{5}\\
&&\nonumber\\
&&(-1)^{|A_1|}A_1(\rho+t)\int^{\rho-\alpha_2+\xi}_{\epsilon_2}[b_{-1},A_2(x+t)]A_3(x-\xi+t)dxA_4(t)\label{6}.
\end{eqnarray}
Now we will rearrange all the terms in such a way that we can express them in terms of known operations.
The terms \rf{5}, \rf{6} can be expressed together as 
\begin{eqnarray}
&&n'_{\rho,\alpha_2, \epsilon_2}(A_1,(A_2,A_3)_{\xi},A_4)(t)+(-1)^{|A_1|}(A_1, (m_{\xi}(A_2, A_3),A_4)_{\epsilon_2})_{\rho}(t)=\nonumber\\
&&n_{\rho,\alpha_2, \epsilon_2}(A_1,(A_2,A_3)_{\xi},A_4)(t)-(m_{\alpha_2}(A_1,(A_2,A_3)_{\xi}),A_4)_{\rho}(t)+\nonumber\\
&&(-1)^{|A_1|}(A_1, (m_{\xi}(A_2, A_3),A_4)_{\epsilon_2})_{\rho}(t).
\end{eqnarray}
Now the second term \rf{2} can be rewritten in terms of $n'$ as follows:
\begin{eqnarray}\label{2sec}
-(-1)^{|A_1||A_2|}n'_{\rho,\alpha_2, \epsilon_1}((A_2,A_1)_{-\alpha_1}, A_3, A_4)(t).
\end{eqnarray}
Lemma 3.1. and Proposition 3.1 allow to represent \rf{2sec} in this way:
\begin{eqnarray}
&&-(-1)^{|A_1||A_2|}n_{\rho,\alpha_2, \epsilon_1}((A_2,A_1)_{-\alpha_1}, A_3, A_4)(t) +\nonumber\\
&&(-1)^{|A_1||A_2|}(m_{\alpha_2}((A_2,A_1)_{-\alpha_1},A_3),A_4)_{\rho}(t)=\nonumber\\
&&\nonumber\\
&&-n_{\rho,\alpha_2,\epsilon_1}((A_1,A_2)_{\alpha_1},A_3,A_4)(t)+\nonumber\\
&&n_{\rho,\alpha_2,\epsilon_1}(Qm_{\alpha_1}(A_1,A_2)+m_{\alpha_1}(QA_1,A_2)(t)+
(-1)^{|A_1|}m_{\alpha_1}(A_1,QA_2), A_3,A_4)(t)+\nonumber\\
&&(-1)^{|A_1||A_2|}(m_{\alpha_2}((A_2,A_1)_{-\alpha_1},A_3),A_4)_{\rho}(t)=\nonumber\\
&&\nonumber\\
&&-n_{\rho,\alpha_2,\epsilon_1}((A_1,A_2)_{\alpha_1},A_3,A_4)+(-1)^{|A_1||A_2|}(m_{\alpha_2}((A_2,A_1)_{-\alpha_1},A_3),A_4)_{\rho}(t)+\nonumber\\
&&((m_{\alpha_1}(A_1,A_2), A_3)_{\alpha_2}, A_4)_{\rho}-(m_{\alpha_1}(A_1,A_2),(A_3,A_4)_{\epsilon_1})_{\rho}(t)-
\nonumber\\
&&Qo^1_{P}(A_1, A_2, A_3, A_4)(t)+o^1_{P}(QA_1, A_2, A_3, A_4)(t)+\nonumber\\
&&(-1)^{|A_1|}o^1_{P}(A_1, QA_2, A_3, A_4)+(-1)^{|A_1|+|A_2|}o^1_{P}(A_1, A_2, QA_3, A_4)(t)+\nonumber\\
&&(-1)^{|A_1|+|A_2|+|A_3|}o^1_{P}(A_1, A_2, A_3, QA_4)(t),
\end{eqnarray}
where $o^1_{P}(A_1,A_2,A_3, A_4)(t)=n_{\rho,\alpha_2,\epsilon_1}(m_{\alpha_1}(A_1,A_2),A_3, A_4)(t)$.

The third term \rf{3} can be expressed as follows:
\begin{eqnarray}\label{3rd}
-(-1)^{|A_1||A_2|+|A_3||A_2|+|A_1||A_3|}(n'_{-\alpha_2,-\xi,-\alpha_1}(A_3,A_2,A_1),A_4)_{\rho}(t).
\end{eqnarray}
The higher commutativity condition, i.e. Lemma 3.2, yields that \rf{3rd} can be rewritten this way:
\begin{eqnarray}
&&(n_{\alpha_2, \alpha_1,\xi}(A_1,A_2,A_3),A_4)_{\rho}(t)-(-1)^{|A_1||A_2|}(m_{\alpha_2}(A_2,A_1)_{-\alpha_1},A_3),A_4)_{\rho}(t)+
\nonumber\\
&&(m_{\alpha_2}(A_1,(A_2,A_3)_{\xi}),A_4)_{\rho}(t)-((m_{\alpha_1}(A_1,A_2),A_3)_{\alpha_2},A_4)_{\rho}-\nonumber\\
&&Qo^2_{P}(A_1, A_2, A_3, A_4)(t)+o^2_{P}(QA_1, A_2, A_3, A_4)(t)+\nonumber\\
&&(-1)^{|A_1|}o^2_{P}(A_1, QA_2, A_3, A_4)(t)+(-1)^{|A_1|+|A_2|}o^2_{P}(A_1, A_2, QA_3, A_4)(t)+\nonumber\\
&&(-1)^{|A_1|+|A_2|+|A_3|}o^2_{P}(A_1, A_2, A_3, QA_4)(t),
\end{eqnarray}
where $o^2_{P}(A_1, A_2, A_3, A_4)(t)=(\t m_{\alpha_2, \alpha_1,\xi}(A_1,A_2,A_3),A_4)_{\rho}$.

Finally, the first and the fourth term together can be written as follows:

\begin{eqnarray}
&&(-1)^{|A_1|}(A_1,n_{\epsilon_2,\xi,\epsilon_1}(A_2,A_3,A_4))_{\rho}(t)-n_{\rho,\alpha_1,\epsilon_2}(A_1,A_2,(A_3,A_4)_{\epsilon_1})(t)-\nonumber\\
&&(-1)^{|A_1|}(A_1,m_{\xi}(A_2,A_3),A_4)_{\epsilon_2})_{\rho}(t)
+(m_{\alpha_1}(A_1,A_2),(A_3,A_4)_{\epsilon_1})_{\rho}(t).
\end{eqnarray}

Summing all up we observe that all the terms containing operation $m$ cancel, and we have the following 
Proposition.\\

\noindent {\bf Proposition 3.3.} {\it Operations $n(\cdot,\cdot,\cdot)$ and $(\cdot,\cdot)$ satisfy the "pentagon relation" in the weak sense:
\begin{eqnarray}
&&Qp_{P}(A_1, A_2, A_3, A_4)(t)-p_{P}(QA_1, A_2, A_3, A_4)(t)-\nonumber\\
&&(-1)^{|A_1|}p_{P}(A_1, QA_2, A_3, A_4)-(-1)^{|A_1|+|A_2|}p_{P}(A_1, A_2, QA_3, A_4)(t)-\nonumber\\
&&(-1)^{|A_1|+|A_2|+|A_3|}p_{P}(A_1, A_2, A_3, QA_4)(t)=\nonumber\\
 &&(-1)^{|A_1|}(A_1,n_{\epsilon_2,\xi,\epsilon_1}(A_2,A_3,A_4))_{\rho}-n_{\rho,\alpha_1,\epsilon_2}(A_1,A_2,(A_3,A_4)_{\epsilon_1})+\nonumber\\
 &&n_{\rho,\alpha_2, \epsilon_2}(A_1,(A_2,A_3)_{\xi},A_4)(t)-n_{\rho,\alpha_2,\epsilon_1}((A_1,A_2)_{\alpha_1},A_3,A_4)(t)+\nonumber\\
&&(n_{\alpha_2, \alpha_1,\xi}(A_1,A_2,A_3),A_4)_{\rho}(t),
\end{eqnarray}
where 
\begin{eqnarray} 
&&p_{P}(A_1, A_2, A_3, A_4)(t)=p'_{P}(A_1, A_2, A_3, A_4)(t)+\nonumber\\
&&n_{\rho,\alpha_2,\epsilon_1}(m_{\alpha_1}(A_1,A_2),A_3, A_4)(t)+(\t m_{\alpha_2, \alpha_1,\xi}(A_1,A_2,A_3),A_4)_{\rho}
\end{eqnarray}
and  the conditions on parameters are: $\rho>>\epsilon_2,\alpha_2>>\epsilon_1,\alpha_1,\xi$.}

\section{Higher order multilinear operations.}
\noindent{\bf 4.1. Short reminder of $A_{\infty}$-algebras.}
The $A_{\infty}$-algebra is a generalization of differential graded associative algebra. Namely, consider a graded vector space $V$ with the differential $Q$. Consider the multilinear operations $\mu_i: V^{\otimes i}\to V$ of the degree $2-i$, such that $\mu_1=Q$. \\

\noindent {\bf Definition 4.1.} (see e.g. \cite{stashbook}){\it The space V is an $A_{\infty}$-algebra if the operations 
$\mu_n$ satisfy bilinear identity:
\begin{eqnarray}\label{arel}
\sum^{n-1}_{i=1}(-1)^{i}M_i\circ M_{n-i+1}=0 
%(-1)^{k+\lambda+k\lambda+nk+k(n_{a_1}+...+n_{a_{\lambda}})}\mu_{n-k+1}
%(a_1, . . . , a_\lambda,
%\mu_k(a_{\lambda+1}, . . . , a_{\lambda+k}), a_{\lambda+k+1}, . . . , a_n) = 0,
\end{eqnarray}
on $V^{\otimes n}$,  
where $M_s$ acts on $V^{\otimes m}$ for any $m\ge s$ as the sum of all possible operators of the form 
${\bf 1}^{\otimes^l}\otimes\mu_s\otimes{\bf 1}^{\otimes^{m-s-l}}$ taken with appropriate signs. In other words, 
\begin{eqnarray}
M_s=\sum^{n-s}_{l=0}(-1)^{l(s+1)}{\bf 1}^{\otimes^l}\otimes\mu_s\otimes{\bf 1}^{\otimes^{m-s-l}}.
\end{eqnarray}}
Let us write several relations which are satisfied by $Q$, $\mu_1$, $\mu_2$, $\mu_3$:
\begin{eqnarray}
&&Q^2=0,\\
&&Q\mu_2(a_1,a_2)=\mu_2(Q a_1,a_2)+(-1)^{|a_1|}\mu_2(a_1,Q a_2),\nonumber\\
&&Q\mu_3(a_1,a_2, a_3)+\mu_3(Q a_1,a_2, a_3)+(-1)^{|a_1|}\mu_3(a_1,Q a_2, a_3)+\nonumber\\
&&(-1)^{|a_1|+|a_2|}\mu_3( a_1, a_2, Q a_3)=\mu_2(\mu_2(a_1,a_2),a_3)-\mu_2(a_1,\mu_2(a_2,a_3)).\nonumber
\end{eqnarray}

In such a way we see that if $\mu_n=0$, $n\ge 3$ , then we have just a differential graded associative algebra (DGA).
If the operations $\mu_n$ vanish for all $n>k$, such $A_{\infty}$-algebras are sometimes called 
$A_{(k)}$-algebras \cite{stasheff}, so e.g. DGA is $A_{(2)}$ algebra.   

We observe that putting $\mu_2\equiv (\cdot, \cdot)$ and $\mu_3={\rm n}$, these relations are manifestly the same as 
the ones relating $Q$, $\mu$ and $\rm n$. Part of the Lian-Zuckerman conjecture is that there are "higher homotopies" $\mu_n, n>3$ satisfying the relations \rf{arel}.   

It is well known that the relations \rf{arel} can be encoded into one equation $\p^2=0$ \cite{stashbook}. To see this one can apply the desuspension operation (the operation which shifts the grading $s^{-1}: V_{n}\to (s^{-1}V)_{n-1}$) to $\mu_n$. In such a way we can define operations of degree $1$: $\t \mu_n=s\mu_n {(s^{-1})}^{\otimes n}$. More explicitly, 
\begin{eqnarray}
\t \mu_n(s^{-1} a_1,...,s^{-1} a_n)=(-1)^{s(a)}s^{-1}\mu_n(a_1,...,a_n),
\end{eqnarray}
such that $s(a)=(1-n)|a_1|+(2-n)|a_2|+...+|a_{n-1}|$. 
The relations between $\t\mu_n$ operations can be summarized in the following simple equations:
\begin{eqnarray}\label{M2}
\sum^n_{i=1}\t M_i\circ \t M_{n+1-i}=0
\end{eqnarray}
on $V^{\otimes n}$, where each $\t M_s$ acts on $V^{\otimes m}$ (for $m\ge s$) as the sum of all operators ${\bf 1}^{\otimes l}\otimes\t \mu_s\otimes{\bf 1}^{\otimes k}$, such that $l+s+k=m$. 
Combining them into one operator $\p=\sum_n\t M_n$, acting on a space $\oplus_kV^{\otimes k}$, the relations \rf{arel} can be summarized in one equation $\p^2=0$.\\ 

\noindent{\bf 4.2. The explicit form of higher operations.} One can check that if we appropriately rename operations from Section 3, namely,  $(\cdot,\cdot)\equiv\mu_2(\cdot,\cdot)$, $n\equiv\mu_3$, $p\equiv\mu_4$ and forget about parameter dependence, we find that the relations between them coincide with $A_{\infty}$ algebra relations. 
We expect that there exist higher order operations related to higher dimensional polytopes, such that they satisfy the $A_{\infty}$-algebra relations modulo parameter dependence. 
We call this object $A_{\infty}$-algebra with parameters. For simplicity in the following we will denote it as $\hat{A}_{\infty}$-algebra.

Here we make a conjecture about the general form of the higher operations of $\hat{A}_{\infty}$-algebra. First of all, we notice that the "commutativity" operations $m$ and $\t m$, which contribute to the multilinear operations, have the following form:
\begin{eqnarray}
&&m_\rho(A,B)(t)=\int^{0}_{-\rho}[b_{-1},(A,B)_{\rho}(t+t')]dt', \nonumber\\
&&\t m_{\rho,\epsilon_1,\epsilon_2}(A,B,C)(t)=\int^{0}_{-\rho}[b_{-1},n'_{\rho,\epsilon_1,\epsilon_2}(A,B,C)](t+t')]dt',
\end{eqnarray}
where $t>\rho$. Therefore, we observe that they depend in a very simple way on the $\hat{A}_{\infty}$-algebra operations of the same multilinear order.
Finally, we notice that the "main" contribution to the 3d and 4th order 
multilinear operations come from $n'$ and $p'$ operations, which correspond to the integral over the first two of Stasheff polytopes (namely, an interval and a pentagon).
Therefore, the general form of the multilinear operation from $\hat{A}_{\infty}$-algebra has the following form:
\begin{eqnarray}
&&\mu^{\rho,K_{n}}_n(A_1,A_2,....,A_n)(t)={\mu'}^{\rho,K_{n}}_n(A_1,A_2,....,A_n)(t),\nonumber\\
&&+\sum_{s}\mu_s^{\epsilon, D_s}(\nu_{n-s}^{\epsilon', D'_{n-s+1}}(A_1, A_2,...,A_s), A_{s+1},...,A_n),
\end{eqnarray}
where by $K_{r}$ we indicate the dependence on the coordinates of a certain $(r-2)$-dimensional Stasheff polytope.  Here $D_s, D'_{k}$ are some Stasheff polytopes of dimensions 
$s-2$, $k-2$ correspondingly, which belong to $\p K_{r}$. The operations $\mu'_n$ and $\nu_n$ are defined as follows:
\begin{eqnarray}
&&{\mu'}^{\rho,K_{n}}_n(A_1,A_2,....,A_n)(t)\equiv\\
&&(-1)^{\frac{(n-3)(n-2)}{2}}(-1)^{(n-2)|A_1|+(n-3)|A_2|+...+|A_{n-2}|}\nonumber\\
&&A_1(\rho+t)\int_{K_{n}}[b_{-1},A_2](t'_2+t)...[b_{-1},A_{n-1}](t'_{n-1}+t)dt'_1\wedge...\wedge dt'_{n-1}A_n(t)\nonumber\\
&&\nu^{\rho,K_{n}}_n(A_1,A_2,....,A_n)(t)=\int^0_{-\rho} [b_{-1},{\mu'}^{\rho,K_{n}}_n(A_1,A_2,....,A_n)(t+t')]dt',\nonumber
\end{eqnarray}
where $K_n$ lies in the domain $\rho>t_2>...>t_n>0$, and 
for $\nu_n$, which is a higher level version of operations $m,\t m$ we have the condition $t>\rho$.
If we act on  ${\mu'}^{\rho,K_{n-2}}_n$ by $Q$, using the Stokes theorem and the properties of boundary of 
Stasheff polytope, we obtain precisely the needed amount of 
terms for the n-th level $A_{\infty}$-algebra relation. 
We hope that the appropriate choice of the polytopes $D_s$, $D'_s$ give us the desired $\hat {A}_{\infty}$-algebra relations. We address this question in the forthcoming article.\\

\noindent {\bf 4.3. Relation to open string theory and real moduli spaces.} In open string theory it is often necessary to  consider the following correlator, see e.g. \cite{dijk}:
\begin{eqnarray}\label{corrs}
\langle \int dt_0d\theta_0\phi(t_0,\theta_0)\int dt_1d\theta_1\phi(t_1,\theta_1)...\int dt_nd\theta_n\phi(t_n,\theta_n)\rangle,
\end{eqnarray}
where $\phi(t,\theta)=\phi^{(0)}(t)+\theta\phi^{(1)}(t)$, $\phi^{(1)}=[b_{-1},\phi^{(0)}]$, where $\phi^{(0)}(t)$ is some open string vertex operator of conformal weight zero and ghost number 1, $\theta$ is an anticommuting operator and the integrals are taken over the supercircle $S^{1|1}$. The correlator \rf{corrs} is obviously not well defined, because it contains the volume of the supergroup generated by $L_{\pm 1}$, $L_0$, $Q$, $b_{\pm 1}$, $b_0$. If we factor it out, we will end up with the integral over the compactified moduli space of $n+1$ points on a circle $\bar{\mathcal{M}}_{n+1}(\mathbb{R})$. However, there is a well known (see e.g. \cite{kap}, \cite{devadoss}) tiling of 
$\bar{\mathcal{M}}_{n+1}(\mathbb{R})$ via Stasheff polytopes. Namely, there are $n!/2$ copies of $K_n$ cells in 
$\bar{\mathcal{M}}_{n+1}(\mathbb{R})$. Consider the following subspace of the integration domain in \rf{corrs}: $t_0=\infty$, $t_1=\rho>0$, $t_n=0$ and $\rho>t_2>t_3>...>t_{n-1}$, where we treat $S^1$ as $\mathbb{R}\cup{\infty}$. It is known that the closure of the inclusion of this space into $\bar{\mathcal{M}}_{n+1}$ is exactly $K_n$. 
Also, we can easily reduce the total integral in \rf{corrs} to this fundamental domain, because the integrand is symmetric with respect to $t_i$.  
Therefore, the natural regularization of \rf{corrs} is the following expression (modulo numerical factor):
\begin{eqnarray}\label{mod}
\lim_{t\to 0}\langle\phi| {\mu'}^{\rho,K_{n}}_n(\phi^{(0)},\phi^{(0)},....,\phi^{(0)})(t)|0\rangle,
\end{eqnarray}
where we used the relation between the polytopes and the multilinear operations $\mu_n'$.  

Hence, $\hat{A}_{\infty}$-algebra relations together with commutativity conditions provide relations between the regularized correlation functions of the type \rf{corrs} in open string theory. It is very important from the point of view of studying the beta-function and conformal perturbation theory. We will return to these questions elsewhere. 

\section{Further remarks}
First of all, we note, that the structures we explored in this article, apply not only to standard TVOAs, but also 
for the conformal field theories with more complicated OPEs, e.g. containing logarithms like in the case of open string. 
This makes the corresponding structures more general than the original homotopy LZ algebra, where the relations strongly relied on the fact that we are dealing with TVOAs only. 

 In this article, we completely neglected the correponding homotopy BV bracket structure, which is defined on TVOA as follows: $\{A,B\}=Res_z[b_{-1},A](z)B$. One immediately can see that this bracket satifies the conditions of the Leibniz algebra, and the symmetrized version of it satisfies the homotopy Lie algebra. As part of the Lian-Zuckerman  conjecture, this homotopy Lie algebra can be extended to $L_{\infty}$-algebra. However, it is impossible to introduce a parameter in this operation as we did with the associative part. 
 What we can do is to take the tensor product of two TVOAs $V$ and $\b V$ and consider the operation which is weakly defined on $V\otimes \b V$ \cite{zeit2}, \cite{zeit3}:
\begin{eqnarray}
\int_{C_{\epsilon,w}}[{\bf b},A(z,\b z)]B(w,\b w),
\end{eqnarray}
where $A,B\in V\otimes \b V$, ${\bf b}=b_{-1}dz+{\b b}_{-1}d\b z$ and $C_{\epsilon,w}$ is a circle of radius $\epsilon$ with the center in $w$. If one expands the corresponding expression in terms of $\epsilon$, at the zero level one obtains the  LZ bracket on the tensor product \cite{azbg}, so that its symmetrized version satisfies the homotopy Jacobi identity.    
Hence, one can expect that the resulting $\epsilon$-dependent bracket will satisfy the homotopy Leibniz/Lie algebra relations modulo parameter dependence. This structures should be related to the geometry of $\bar{M}_{0,n}$ moduli space and constructions of string field theory \cite{zwiebach}.

Finally, we note, that recently it was observed that the Maurer-Cartan equations for the Lian-Zuckerman homotopy algebras associated with semi-infinite complex of Virasoro algebra lead to classical field equations and their string corrections \cite{lmz},\cite{bvym}, \cite{zeit2}, \cite{zeit3}, \cite{cftym}, \cite {khromov}. It will be extremely important to define the proper anolgue of the Maurer-Cartan equation and its symmetries in the parameter-dependent case. 

\section{Acknowledgements}
I am very grateful to I.B. Frenkel, M.M. Kapranov and R. Raj for useful discussions on the subject. I am indebted to Jim Stasheff for reading the manuscript and for many useful comments and suggestions.

\end{document}